\documentclass[12pt]{article}
\usepackage[utf8]{inputenc}
\usepackage{amssymb, amsmath,amsthm, thmtools, mathabx, mathtools}
\usepackage{bbm}
\usepackage{amsfonts}
\usepackage[left=2.5cm,right=2.5cm,top=2cm,bottom=2cm]{geometry}
\usepackage{xcolor}

\usepackage{tabto}
\TabPositions{2 cm, 4 cm, 6 cm, 8 cm}

\usepackage{tikz-cd}

\usepackage{tikz}
\usetikzlibrary{arrows, automata, positioning}

\usepackage{tocloft}
\usepackage{appendix}

\usepackage{enumerate}

\usepackage{hyperref}
\usepackage[capitalize]{cleveref}

\hypersetup{
  colorlinks   = true, 
  urlcolor     = blue, 
  linkcolor    = blue, 
  citecolor   = blue 
}

\newtheorem{theorem}{Theorem}[section]

\newtheorem{proposition}[theorem]{Proposition}

\newtheorem*{theorem*}{Theorem}

\theoremstyle{definition}

\newtheorem{remark}[theorem]{Remark}

\newtheorem*{remark*}{Remark}
\newtheorem*{notation*}{Notation}
\newtheorem*{acks*}{Acknowledgements}
\newtheorem*{out*}{Outline}

\renewcommand\leq{\leqslant}

\newcommand{\T}{\mathrm{(T)}}

\begin{document}

\title{A torsion-free non-sofic group}
\author{Francesco Fournier-Facio}
\date{\today}
\maketitle

\begin{abstract}
OpenAI announced the existence of a non-sofic group: the unit group of the binary Leavitt algebra. We exhibit a different source of examples (relying on the same technical criterion) which includes torsion-free groups.
\end{abstract}

\section{Introduction}

On 1 August 2026, OpenAI announced solutions to 10 major open problems in mathematics. In group theory, the target was the problem of the existence of a non-sofic group, attributed to Gromov \cite{gromov} and Weiss \cite{weiss}. We refer the reader to \cite[Chapter 7]{CSC} for the relevant definitions.

\begin{theorem}[\cite{openai}]
There exists a non-sofic group.
\end{theorem}

The key technical ingredient in the proof is the following, which crucially builds on work of Kun \cite{kun} and Kun--Thom \cite{kunthom}\footnote{Contrary to the claim in the OpenAI announcement on 1 August \cite{openai} that there has been ``no progress on the main result for at least a decade''. This claim has been redacted on 3 August.}.

\begin{proposition}[{\cite[Proposition 2.3]{openai}}]
\label{criterion}
Let $\Gamma \leq G$ be infinite groups with property $\T$. Suppose that $G = \langle \Gamma, t_1, \ldots, t_m \rangle$ where $t_i \Gamma t_i^{-1} \subset \Gamma$ for all $i$. Suppose that there exists a finitely generated group $J \leq G$ such that $[\Gamma, J] = \Gamma \cap J = \{ 1 \}$ and $t_1 J t_1^{-1} \subset \Gamma$. If $G$ is sofic, then $J$ is LEF.
\end{proposition}

This proposition is then applied with $J$ being Thompson's group $V$ (which is not LEF, being finitely presented and simple), and $G$ and $\Gamma$ being groups of elementary matrices over the binary Leavitt algebra $R$. This is a clever choice, because the behaviour in the proposition is very natural for Thompson-like groups, but Thompson-like groups almost never have property $\T$. On the other hand, elementary matrices over $R$ do, by work of Ershov--Jaikin-Zapirain \cite{ershovjaikin}, because $R$ is a finitely generated ring, and they can be seen as living in the same world as $V$, acting on binary strings. To the best of my knowledge, this connection between Leavitt algebras and Thompson groups was first observed independently by Birget \cite{birget} and Nekrashevych \cite{nekr}.

\medskip

In an effort to better understand the solution, and the role that the various parts of the proof play, I explain here how to construct a non-sofic group avoiding Leavitt algebras entirely, while still using Proposition \ref{criterion}. I do not claim that this proof is easier, better, or more elementary, although it is more natural to me, in that it only uses things that I already knew on 30 July (unlike the connection between Leavitt algebras and Thompson groups). It is a rather standard argument in small cancellation theory \cite{osin, hull}, which to me indicates that Proposition \ref{criterion} is the key novelty in \cite{openai}. As it often happens with these things, the novelty is not so much in the proof (which is still clever and technical, as beautifully explained by Andreas Thom \cite{andreas}), as in the statement: experts aware of \cite{kunthom} would probably have been able to prove Proposition \ref{criterion}, but they might not have thought of this statement. There are myriad possible statements strengthening \cite{kunthom}, many of which are surely false, and many of which might be true but have no application. AI models are not inconvenienced by constraints of time and resources, so they get to throw things at the wall and see what sticks.

\medskip

In the next section I present this alternative construction, which buys some new flexibility:

\begin{theorem}
There exists a finitely presented torsion-free non-sofic group.
\end{theorem}

It is easy to see that for a Leavitt algebra $R$ like the one used in \cite{openai}, the multiplicative group $R^\times$ (hence any group of elementary matrices over $R$ \cite[Corollary 4.4]{matrixgen}) always has torsion. The existence of a non-sofic group immediately implies the existence of a finitely presented one, because non-soficity is an open property in the space of marked groups \cite[Proposition 7.5.10]{CSC}, and so \cite{openai} already yields a finitely presented non-sofic group. But this argument does not specialise to torsion-free groups, since torsion-freeness is not an open property; here we directly construct a finitely presented example.

\begin{acks*}
I am supported by the Herchel-Smith Postdoctoral Fellowship Fund.

On 31 July 2026 I saw an early version of \cite{openai} and went through it together with Henry Bradford and Alon Dogon. I posted the first version of this note on my homepage on the evening of 1 August 2026, the day of the announcement. Since then, I received several comments that improved it significantly: Goulnara Arzhantseva and Marcin Kotowski pointed to missing citations, Johan {\"O}inert explained that $R^\times$ has torsion even if the base field is infinite, him and Zixiang Zhou noticed a mathematical typo in the previous paragraph, Benjamin Steinberg taught me about the first instances of the Leavitt--Thompson connection, and Matt Zaremsky simplified the last step of the proof. This is an essential part of the process of doing mathematics, and it is possible only when the author is named and lists an email address whereby anyone can contact them\footnote{Compare with \cite{openai}.}.
\end{acks*}

\begin{remark}
On 7 August 2026, a new Kun--Thom paper appeared on the arXiv \cite{kunthom2}, which builds on \cite[Proposition 2.3]{openai} to give a more robust construction of non-sofic groups, in the form of wreath products. This also shows that there are weakly sofic groups that are not sofic \cite{glebsky}, clarifying the relationship between soficity and the other approximation properties for which non-examples are still to be found. The proof relies on a centraliser rigidity theorem by Alekseev--Thom \cite{alekseevthom} (further developing the ideas from \cite{kunthom}, in the context of groupoids) which also appeared on 7 August 2026 on the arXiv, but which had been circulating for a month already. We leave the rest of this note as it was before this new development.
\end{remark}

\pagebreak

\section{The construction}

We will use the following input groups:
\begin{itemize}
\item A universal finitely presented torsion-free group $U$, that is, a finitely presented torsion-free group that contains a copy of every finitely presented torsion-free group \cite{higman, unitf1, unitf2};
\item A finitely presented simple torsion-free group $S$, e.g. a Burger--Mozes \cite{BM} or Hyde--Lodha group \cite{HL};
\item A torsion-free hyperbolic property $\T$ group $H$, e.g. a random group at suitable density \cite{zuk, kotowski}.
\end{itemize}

The group $U$ embeds into a finitely presented torsion-free property $\T$ group $P$. This is essentially a consequence of SQ-universality of hyperbolic groups \cite{gromov:hyp, delzant, olsh} and it can be achieved via a small cancellation theory argument, as in \cite{AMO}: consider the relatively hyperbolic pair $(U * H, U)$, and add relations identifying each generator of $U$ with an element of $H$; see \cite[Proposition 2.3]{approx} for a precise statement. This step can be done while preserving torsion-freeness, see \cite[Theorem 2.4.5]{osin}.

\medskip

By universality, $P$ contains a subgroup of the form $P_1 \times P_2 \times S \leq P$, where each $P_i$ is isomorphic to $P$. Let $E$ denote the double HNN extension of $P$ with stable letters $u_1, u_2$, where $u_i$ conjugates $P$ to $P_i$. Note that $E$ is torsion-free.

\medskip

Consider the action of $E$ on its Bass--Serre tree. Since $P_1$ and $P_2$ are disjoint edge groups, there is a path of length $2$ with trivial stabiliser. By \cite{MO}, the group $E$ is acylindrically hyperbolic. We now perform one more small cancellation step: there exists a finitely presented torsion-free quotient $\pi \colon E \to G$ such that $G$ has property $\T$, and $\pi(S) \neq \{ 1 \}$: this is an application of \cite[Corollary 7.4]{hull}, which gives a common quotient of $E$ and $H$ which is injective on a given finite set. Also here one can preserve torsion-freeness, see \cite[Theorem 7.1(e)]{hull}, and finite presentability: since both $E$ and $H$ are finitely generated, only finitely many relators need to be added to obtain a quotient of $E * H$ onto which both $E$ and $H$ surject. Since $S$ is simple, it follows that $\pi|_S$ is injective.

\medskip

Let $\Gamma = \pi(P) \leq G$ and $t_i = \pi(u_i)$. Being a quotient of $P$, $\Gamma$ also has property $\T$, and it is infinite because it contains $\pi(S) \cong S$. Moreover $G = \langle \Gamma, t_1, t_2 \rangle$ and $t_i \Gamma t_i^{-1} \subset \Gamma$, because the analogous statements hold for $E = \langle P, u_1, u_2 \rangle$. Let $J = t_1^{-1} \pi(S) t_1$. Since $[P_1, S] = 1$, we have $[\Gamma, J] = 1$. Finally, we claim that $\Gamma \cap J = \{ 1 \}$: since $\Gamma$ and $J$ commute, an element in the intersection belongs to the centre of $J$, which is trivial because $J \cong S$ is simple.

\medskip

Therefore $G = \langle \Gamma, t_1, t_2 \rangle$ and $J$ satisfy the hypotheses of Proposition \ref{criterion}. If $G$ is sofic, then $J$ is LEF, but this is impossible since $J \cong S$ is finitely presented, infinite and simple, and finitely presented LEF groups are residually finite \cite[Proposition 7.3.8]{CSC}.

\pagebreak

\footnotesize

\bibliographystyle{amsalpha}
\bibliography{ref}

\providecommand{\bysame}{\leavevmode\hbox to3em{\hrulefill}\thinspace}
\providecommand{\MR}{\relax\ifhmode\unskip\space\fi MR }
\providecommand{\MRhref}[2]{%
  \href{http://www.ams.org/mathscinet-getitem?mr=#1}{#2}
}
\providecommand{\href}[2]{#2}
\begin{thebibliography}{AMO07}

\bibitem[AMO07]{AMO}
G.~Arzhantseva, A.~Minasyan, and D.~Osin, \emph{The {SQ}-universality and
  residual properties of relatively hyperbolic groups}, J. Algebra \textbf{315}
  (2007), no.~1, 165--177. \MR{2344339}

\bibitem[AT26]{alekseevthom}
V.~Alekseev and A.~Thom, \emph{Centralizers of sofic approximations of kazhdan
  groups}, arXiv preprint arXiv:2608.05362, 2026.

\bibitem[Bir04]{birget}
J.-C. Birget, \emph{The groups of {R}ichard {T}hompson and complexity},
  vol.~14, 2004, International Conference on Semigroups and Groups in honor of
  the 65th birthday of Prof. John Rhodes, pp.~569--626. \MR{2104771}

\bibitem[BM97]{BM}
M.~Burger and S.~Mozes, \emph{Finitely presented simple groups and products of
  trees}, C. R. Acad. Sci. Paris S\'{e}r. I Math. \textbf{324} (1997), no.~7,
  747--752. \MR{1446574}

\bibitem[BS08]{unitf1}
I.~Belegradek and A.~Szczepa\'{n}ski, \emph{Endomorphisms of relatively
  hyperbolic groups}, Internat. J. Algebra Comput. \textbf{18} (2008), no.~1,
  97--110, With an appendix by O. V. Belegradek. \MR{2394723}

\bibitem[Chi14]{unitf2}
M.~Chiodo, \emph{On torsion in finitely presented groups}, Groups Complex.
  Cryptol. \textbf{6} (2014), no.~1, 1--8. \MR{3200358}

\bibitem[CSC23]{CSC}
T.~Ceccherini-Silberstein and M.~Coornaert, \emph{Cellular automata and
  groups}, Springer Monographs in Mathematics, Springer, Cham, [2023]
  \copyright 2023, Second edition [of 2683112]. \MR{4738386}

\bibitem[Del96]{delzant}
T.~Delzant, \emph{Sous-groupes distingu\'{e}s et quotients des groupes
  hyperboliques}, Duke Math. J. \textbf{83} (1996), no.~3, 661--682.
  \MR{1390660}

\bibitem[EJZ10]{ershovjaikin}
M.~Ershov and A.~Jaikin-Zapirain, \emph{Property ({T}) for noncommutative
  universal lattices}, Invent. Math. \textbf{179} (2010), no.~2, 303--347.
  \MR{2570119}

\bibitem[FF25]{approx}
F.~Fournier-Facio, \emph{Stability, approximable quotients, and higher property
  {(T)}}, arXiv preprint arXiv:2512.09180, 2025.

\bibitem[Gle23]{glebsky}
L.~Glebsky, \emph{Extensions of a residually finite group by a weakly sofic
  group are weakly sofic}, Rev. Mat. Iberoam. \textbf{39} (2023), no.~3,
  1097--1104. \MR{4603646}

\bibitem[Gro87]{gromov:hyp}
M.~Gromov, \emph{Hyperbolic groups}, Essays in group theory, Math. Sci. Res.
  Inst. Publ., vol.~8, Springer, New York, 1987, pp.~75--263. \MR{919829}

\bibitem[Gro99]{gromov}
\bysame, \emph{Endomorphisms of symbolic algebraic varieties}, J. Eur. Math.
  Soc. (JEMS) \textbf{1} (1999), no.~2, 109--197. \MR{1694588}

\bibitem[Hig61]{higman}
G.~Higman, \emph{Subgroups of finitely presented groups}, Proc. Roy. Soc.
  London Ser. A \textbf{262} (1961), 455--475. \MR{130286}

\bibitem[HL25]{HL}
J.~Hyde and Y.~Lodha, \emph{Finitely presented simple left-orderable groups in
  the landscape of {R}ichard {T}hompson's groups}, Ann. Sci. \'{E}c. Norm.
  Sup\'{e}r. (4) \textbf{58} (2025), no.~2, 419--432. \MR{4902401}

\bibitem[Hul16]{hull}
M.~Hull, \emph{Small cancellation in acylindrically hyperbolic groups}, Groups
  Geom. Dyn. \textbf{10} (2016), no.~4, 1077--1119. \MR{3605028}

\bibitem[KK13]{kotowski}
M.~Kotowski and M.~Kotowski, \emph{Random groups and property {$(T)$}:
  \.{Z}uk's theorem revisited}, J. Lond. Math. Soc. (2) \textbf{88} (2013),
  no.~2, 396--416. \MR{3106728}

\bibitem[KT19]{kunthom}
G.~Kun and A.~Thom, \emph{Inapproximability of actions and {K}azhdan's property
  {(T)}}, arXiv preprint arXiv:1901.03963, 2019.

\bibitem[KT26a]{matrixgen}
H.~V. Khanh and V.~H. Thanh, \emph{Matrix generators for the unit groups of
  $l_k (1, d) $}, arXiv preprint arXiv:2607.10351, 2026.

\bibitem[KT26b]{kunthom2}
G.~Kun and A.~Thom, \emph{Nonsofic wreath products of residually finite
  groups}, arXiv preprint arXiv:2608.06222, 2026.

\bibitem[Kun16]{kun}
G.~Kun, \emph{On sofic approximations of property {(T)} groups}, arXiv preprint
  arXiv:1606.04471, 2016.

\bibitem[MO15]{MO}
A.~Minasyan and D.~Osin, \emph{Acylindrical hyperbolicity of groups acting on
  trees}, Math. Ann. \textbf{362} (2015), no.~3-4, 1055--1105. \MR{3368093}

\bibitem[Nek04]{nekr}
V.~V. Nekrashevych, \emph{Cuntz-{P}imsner algebras of group actions}, J.
  Operator Theory \textbf{52} (2004), no.~2, 223--249. \MR{2119267}

\bibitem[Ols93]{olsh}
A.~Yu. Olshanskii, \emph{On residualing homomorphisms and {$G$}-subgroups of
  hyperbolic groups}, Internat. J. Algebra Comput. \textbf{3} (1993), no.~4,
  365--409. \MR{1250244}

\bibitem[Ope26]{openai}
OpenAI, \emph{Non-sofic groups exist},
  \url{https://openai.com/index/ten-advances-in-mathematics/}, 2026.

\bibitem[Osi10]{osin}
D.~Osin, \emph{Small cancellations over relatively hyperbolic groups and
  embedding theorems}, Ann. of Math. (2) \textbf{172} (2010), no.~1, 1--39.
  \MR{2680416}

\bibitem[Tho]{andreas}
A.~Thom, \emph{What are the key new ideas in the proof of nonsoficity of groups
  in {OpenAI}’s construction of nonsofic groups}, MathOverflow,
  \url{https://mathoverflow.net/q/513885}.

\bibitem[Wei00]{weiss}
B.~Weiss, \emph{Sofic groups and dynamical systems}, vol.~62, 2000, Ergodic
  theory and harmonic analysis (Mumbai, 1999), pp.~350--359. \MR{1803462}

\bibitem[{\.Z}uk03]{zuk}
A.~{\.Z}uk, \emph{Property ({T}) and {K}azhdan constants for discrete groups},
  Geom. Funct. Anal. \textbf{13} (2003), no.~3, 643--670. \MR{1995802}

\end{thebibliography}

\vspace{0.5cm}

\normalsize

\noindent{\textsc{Department of Pure Mathematics and Mathematical Statistics, University of Cambridge, UK}}

\noindent{\textit{E-mail address:} \texttt{ff373@cam.ac.uk}}

\end{document}